\documentclass[english,11pt]{mysmfart}

\usepackage[utf8]{inputenc}
\usepackage[a4paper,margin=3cm]{geometry}
\usepackage{amssymb,footnote}
\usepackage{multicol}
\usepackage{multirow}
\usepackage{dsfont}
\usepackage{bm}
\usepackage{comment}
\usepackage{microtype}
\usepackage[colorlinks=true,linkcolor=black,citecolor=black,hypertexnames=false]{hyperref}
\usepackage{caption}
\usepackage[ruled,vlined,boxed]{algorithm2e}
\usepackage{stmaryrd}
\usepackage{color}
\usepackage[all]{xy}
\usepackage{enumerate}
\usepackage{tikz}
\usepackage{url}
\usepackage[english]{babel}

\newenvironment{myenumerate}[1][1.]{%
  \vspace{-\parskip}
  \begin{enumerate}[#1]\setlength{\itemsep}{0.2ex}}
 {\end{enumerate}%
  \vspace{-\parskip}}
\newenvironment{myitemize}{%
  \vspace{-\parskip}
  \begin{itemize}\setlength{\itemsep}{0.2ex}}
 {\end{itemize}%
  \vspace{-\parskip}}

\newcommand{\ph}{$\vphantom{A^A_A}$}

\newtheorem{prop'}{Proposition}[section]
\newtheorem{conj'}[prop']{Conjecture}
\newtheorem{lem'}[prop']{Lemma}
\newtheorem{thm'}[prop']{Theorem}
\newtheorem{cor'}[prop']{Corollary}
\theoremstyle{definition}
\newtheorem{hyp'}[prop']{Hypotheses}
\newtheorem{definit'}[prop']{Definition}
\newtheorem{ex'}[prop']{Example}
\newtheorem{ques'}[prop']{Question}
\newtheorem{rem'}[prop']{Remark}

\theoremstyle{plain}
\newtheorem{prop}{Proposition}
\newtheorem{conj}[prop]{Conjecture}

\theoremstyle{definition}

\def\={\buildrel {\rm d\acute ef}\over =}

\DeclareMathOperator{\Card}{Card}

\DeclareMathOperator{\Spec}{Spec}

\DeclareMathOperator{\Aut}{Aut}
\DeclareMathOperator{\Ind}{Ind}
\DeclareMathOperator{\Hom}{Hom}
\DeclareMathOperator{\Gal}{Gal}
\newcommand{\GGal}{\text{\bf Gal}}

\DeclareMathOperator{\Sym}{Sym}

\DeclareMathOperator{\GL}{GL}
\newcommand{\nr}{\text{\rm nr}}

\newcommand{\norm}{\vert\cdot\vert}

\renewcommand{\P}{\mathbb{P}}
\newcommand{\GG}{\mathcal{G}}

\newcommand{\RR}{\mathcal{R}}

\newcommand{\JJ}{\mathfrak{J}}

\newcommand{\MK}{\mathfrak{M}}

\newcommand{\bX}{\mathbb{X}}

\newcommand{\oK}{{\mathcal O}_K}

\renewcommand{\1}{{\mathbb 1}}
\newcommand{\Z}{{\mathbb Z}}
\newcommand{\N}{{\mathbb N}}
\newcommand{\Q}{{\mathbb Q}}

\newcommand{\Qp}{\Q_{p}}

\newcommand{\Zp}{\Z_{p}}

\newcommand{\ur}{\text{\rm ur}}
\newcommand{\tr}{\text{\rm tr}}

\newcommand{\F}{\mathbb F}
\newcommand{\Fp}{{\mathbb F}_{p}}

\newcommand{\R}{\mathbb R}

\newcommand{\Qbar}{\overline\Q}
\newcommand{\Qpbar}{\Qbar_p}
\newcommand{\Zpbar}{\overline \Z_p}
\newcommand{\Fpbar}{\overline \F_p}
\newcommand{\rhobar}{\overline\rho}

\newcommand{\D}{\mathcal{D}}

\newcommand{\ttt}{{\rm t}}

\def\DD{{\mathcal D}}

\def\WW{{\mathcal W}}

\newcommand{\Res}{\text{\rm Res}}

\newcommand{\aRes}{\text{\rm aRes}}
\newcommand{\mRes}{\text{\rm mRes}}
\newcommand{\Sd}{\mathfrak S_d}

\newcommand{\Ga}{\mathbb G_a}
\newcommand{\Gm}{\mathbb G_m}

\newcommand{\Vect}[1]{#1\textbf{-Vect}}
\newcommand{\Alg}[1]{#1\textbf{-Alg}}
\newcommand{\Sch}[1]{#1\textbf{-Sch}}

\newcommand{\GR}{\GG\RR}
\newcommand{\GRbar}{\overline{\GR}}

\newcommand{\gA}{\text{\tt A}}
\newcommand{\gB}{\text{\tt B}}
\newcommand{\gAB}{\text{\tt AB}}
\newcommand{\gO}{\text{\tt O}}

\newcommand{\I}{\text{\rm I}}
\newcommand{\II}{\text{\rm II}}

\renewcommand{\epsilon}{\varepsilon}

\author[X. Caruso]{Xavier Caruso}
\address{CNRS; IMB,
Université de Bordeaux,
351 cours de la Libération,
33405 Talence, France}
\email{xavier.caruso@normalesup.org}

\author[A. David]{Agnès David}
\address{LMB,
Universit\'e de Franche-Comt\'e,
16 route de Gray,
25030 Besançon Cedex, France;
IRMAR,
Université de Rennes I,
Campus de Beaulieu,
35042 Rennes Cedex, France}
\email{agnes.david@math.cnrs.fr}

\author[A. Mézard]{Ariane Mézard}
\address{DMA,
École Normale Supérieure PSL,
45 rue d'Ulm,
75005 Paris, France}
\email{ariane.mezard@ens.fr}

\title[Can we dream of a 1-adic Langlands correspondence?]
{Can we dream of a\\1-adic Langlands correspondence?}

\begin{document}

\begin{abstract}
After observing that some constructions and results in the
$p$-adic Langlands programme are somehow independent from~$p$,
we formulate the hypothesis that this astonishing uniformity 
could be explained by a $1$-adic Langlands correspondence.
\end{abstract}

\maketitle

\hfill{\it To Catriona Byrne.}

\tableofcontents

The Langlands programme is a far-reaching and influential web of 
theorems and conjectures which has motivated a lot of research in 
Number Theory and Arithmetic Geometry for more than fifty years. Very 
roughly, it stipulates a profound and meaningful correspondence between 
representations of Galois groups on the one hand and representations of 
reductive groups on the other hand.
Many variations on this theme are actually possible, depending on which 
base field (number field, function field, $p$-adic field, \emph{etc}.)
and which category of coefficients we are working with.

The pioneer works in the Langlands programme were mostly concerned
with $\mathbb C$-valued representations. However,
since the beginning of the 21st century, a purely $p$-adic version of 
Langlands correspondence has emerged under the impulsion of Breuil.
Nowadays, this $p$-adic correspondence is fully established for 
$2$-dimensional representations of 
$\Gal(\Qpbar/\Qp$) but, beyond this, little is known. 
Many examples have however been worked out and several conjectures 
have been proposed---and sometimes proved---throughout the years. 
One of them is the Breuil-Mézard conjecture, which predicts that 
the geometrical properties of some Galois deformations spaces 
are directly related to the decomposition properties of some
representations of the corresponding reductive group.

Looking more carefully into the aforementioned works, we 
notice that, in many cases, the underlying prime number~$p$ often
plays a figurative role in the calculations.
Typically, the relevant reductive groups are usually defined over
$\Z$ and a significant part of the constructions and arguments can
be carried out at this level. On the Galois side, this constancy
is not so obvious but it is nevertheless visible; indeed, even 
though a prime number needs to be fixed from the very beginning, 
we often observe, at the end of the day, that the results of the
computations are mostly independent from it.

We make the hypothesis that these strong properties 
of uniformity with respect to~$p$ could have a deep meaning and all 
be the consequences of a new type of Langlands correspondence, which 
should be considered as the common denominator of the $p$-adic 
Langlands correspondences when $p$ varies.
We call this new hypothetic correspondence the \emph{$1$-adic 
Langlands correspondence} because we believe that the natural language 
to formulate it is the mysterious theory of characteristic one whose 
main protagonist is the famous field with one element.

The aim of this note is to bring the reader to the agreement that 
our hypothesis is not crazy but has conceivable foundations and
deserves consideration. We start our argumentation 
by reviewing in \S \ref{sec:BM} some recent developments towards the 
Breuil-Mézard conjecture, with the objective to highlight the places 
where the arguments and/or the notions take a combinatorial flavour 
in which we have the feeling that the underlying prime number~$p$ 
plays a secondary role.
Then, in \S \ref{sec:f1}, we briefly recall the philosophy of the
field with one element and show that it is incredibly appropriate 
for interpreting many objects and carrying out many constructions
encountered in \S \ref{sec:BM}.
Finally, in order to give more substance to our dream, we conclude 
this article by an appendix in which we share some thoughts towards
the development of a Galois theory in characteristic one, which is
certainly a prerequisite for a $1$-adic Langlands correspondence.

\section{Combinatorics around the Breuil-Mézard conjecture}
\label{sec:BM}

Let $p>2$ be a prime number.
Throughout this section, we fix a finite extension $K$ of $\Qp$
and write $G_K=\Gal(\overline{\Q}_p/K)$ for its absolute Galois
group.
The Breuil-M\'ezard conjecture is a concrete statement relating 
deformations spaces of representations of $G_K$, on the one hand,
and representations of $p$-adic reductive groups, on the other
hand. The aim of this section is, firstly, to recall the 
formulation of this conjecture and, secondly, to emphasize that,
in many cases, it can be approched using combinatorial arguments 
and constructions.
These observations will be the key to build bridges with the
field with one element in \S \ref{sec:f1}.

\subsection{Review on the Breuil-Mézard conjecture}

We denote by $\oK$ (resp $k_K$) the ring of integers (resp. the
residue field) of $K$.
Let $\rhobar : G_K \to \GL_n(\Fpbar)$ be a continuous 
$\Fpbar$-representation 
of $G_K$ of dimension $n$ and let $R_{\rhobar}$ denote the
$\Zpbar$-algebra parametrizing the deformations of $\rhobar$. 
In~\cite{Ki3}, Kisin proved that $R_{\rhobar}$ admits quotients
with strong arithmetical interest. More precisely, given in
addition the two following data:
\begin{myitemize}
\item a \emph{Hodge type} $\lambda$, that is, by definition, the
datum of a tuple $(\lambda_1, \ldots, \lambda_n) \in \Z^n$ with 
$\lambda_1 \geq \cdots \geq \lambda_n$ for all embedding $\iota:
k_K \hookrightarrow \Fpbar$,
\item an \emph{inertial type} $\ttt$, that is, by definition, a 
finite dimensional $\Qpbar$-representation of the inertia
subgroup $I_K \subset G_K$ having open kernel and admitting an 
extension to $G_K$,
\end{myitemize}
Kisin constructed a surjective morphism of $\Zpbar$-algebras 
$R_{\rhobar} \to R^{\lambda,\ttt}_{\rhobar}$ that parametrizes the 
lifts of $\rhobar$ which are potentially crystalline with Hodge Tate 
weights $(\lambda_i + n - i)_{i,\, \iota}$ and inertial type $\ttt$.

The Breuil-Mézard conjecture is a numerical relation between the 
Hilbert-Samuel multiplicity of the special fibre of
$R^{\lambda,\ttt}_{\rhobar}$, denoted by $e(R^{\lambda,\ttt}_{\rhobar}
\otimes_{\Zpbar} \Fpbar)$, and
invariants coming from the representation theory of $\GL_n$.
Precisely, let $L_\lambda$ be the $\Zpbar$-representation of 
$\GL_n(\oK)$ of highest weight $\lambda$. After \cite{He,SZ,CEG+},
we know that there is a finite dimensional smooth irreducible 
$\Qpbar$-representation $\sigma(\ttt)$ of $\GL_n(\oK)$ associated
to $\ttt$. We choose a $\GL_n(\oK)$-stable $\Zpbar$-lattice 
$L_{\ttt}$ in $\sigma(\ttt)$, form the tensor product
$L_{\lambda,\ttt} = L_{\ttt}\otimes_{\Zpbar} L_{\lambda}$ and write
its semi-simplification modulo $p$ as follows:
$$\big(L_{\lambda,\ttt}\otimes_{\Zpbar} \Fpbar\big)^{\text{ss}} 
\simeq  \bigoplus_{\sigma\in\DD} \sigma^{n_{\lambda,\ttt}(\sigma)}$$
where the sums runs over the set $\DD$ of Serre weights $\sigma$, that 
is the set of (isomorphism classes of) irreducible $\Fpbar$-representations of $\GL_n(k_K)$.

\begin{conj}[Breuil-Mézard]
\label{conjBM}
There exists a family of integers $(\mu_{\rhobar}(\sigma))_{\rhobar,
\sigma}$, called \emph{intrinsic multiplicities}, such that the 
following numerical equality holds:
\begin{equation}
\label{eq:BM}
e\big(R_{\rhobar}^{\lambda,\ttt} \otimes_{\Zpbar} \Fpbar\big) =
\sum_{\sigma\in\DD} n_{\lambda,\ttt}(\sigma)\:\mu_{\rhobar}(\sigma)
\end{equation}
for all triple $(\rhobar,\lambda,\ttt)$ as above.
\end{conj}

The Breuil-Mézard conjecture was first formulated for
$2$-dimensional representations in~\cite{BM1}. Since then, it has
attracted a lot of attention. Kisin~\cite{Ki1} proved it when $K =
\Qp$ (and $n = 2$) by making use of the $p$-adic local Langlands 
correspondence for $\GL_2(\Qp)$ and the (global) Taylor-Wiles-Kisin 
patching argument. Sander~\cite{Sa} and Paškūnas~\cite{Pa} gave a 
purely local alternative proof which has been extended later on
by Hu and Tan to nonscalar split residual representations~\cite{HuTa}.
For a general $K$, but still assuming $n = 2$, the conjecture
was proved by Gee and Kisin~\cite{GeKi} (see also 
\cite[Appendix~C]{CEGS}) when $\lambda = (0, 0)$ for each
embedding (which corresponds to potentially Barsotti--Tate
deformations).

The extension of the Breuil-Mézard conjecture to higher~$n$ came later. 
The case of $3$-dimensional representations was considered and 
partially solved by Herzig, Le and Morra~\cite{HLM} and Le, Le Hung, 
Levin and Morra~\cite{LLHLM1,LLHLM}. The formulation stated in dimension $n$ in
Conjecture~\ref{conjBM} is due to Emerton and Gee~\cite{EmGe}
(see also~\cite{GHS}).

\subsection{Encoding representations with combinatorial data}
\label{ssec:encode}

It turns out that, in some cases, all the objects that intervene in 
the statement of the Breuil-Mézard conjecture can be entirely described
by combinatorial data. Besides, these explicit descriptions provide us
with a new viewpoint on the conjecture, quite useful for
attacking it.

\paragraph*{The case of $2$-dimensional representations}

We start with the case of $\GL_2$ for which encodings are simpler
and more is known. For simplicity, we assume further that $K/\Qp$ is
unramified, \emph{i.e.} $K = \Q_{p^f}$ for some positive integer~$f$.
In dimension $2$, the \emph{irreducible} continuous representations
$\rhobar : G_K \to \GL_2(\Fpbar)$ all take the form:
\begin{equation}
\label{eq:rhobar}
\rhobar =
\Ind_{G_{K'}}^{G_K} \Big( \omega_{2f}^h \cdot \nr'(\theta)\Big)
\end{equation}
where $K'$ is the unique unramified extension of degree $2$ of $K$,
$\omega_{2f}$ is the fundamental character of $G_{K'}$ of level $2f$
and $\nr'(\theta)$ denotes the unique unramified character of $G_{K'}$ 
sending the arithmetic Frobenius to $\theta$. The parameters $h$ and
$\theta$ are an integer defined modulo $p^{2f}{-}1$ and an element
of $\Fpbar$ respectively.

Similarly, we have a complete description of \emph{tame} inertial 
types, that are, by definition, those inertial types $\ttt : I_K \to 
\GL_2(\Fpbar)$ that factor through the tame inertia. Depending on
whether they are reducible or not, they are of the form:
\begin{equation}
\label{eq:ttt}
\ttt = \omega_f^\gamma \oplus \omega_f^{\gamma'}
\quad \text{or} \quad
\ttt = \Ind_{I_{K'}}^{I_K} \omega_{2f}^\gamma
\end{equation}
with obvious notations. In the first case, we say that $\ttt$ has
level~$f$; otherwise, that it has level~$2f$.

We now assume further that $\lambda = (0,0)$ for each embedding.
In this case, the integers $n_{\lambda,\ttt}(\sigma)$ are always
$0$ or $1$ and we define $\DD(\ttt)$ as the set of Serre weights
$\sigma \in \DD$ for which $n_{\lambda,\ttt}(\sigma) = 1$.
Similarly, it is conjectured that $\mu_{\rhobar}(\sigma)$ is
always $0$ or $1$ as well and we let $\DD(\rhobar) \subset \DD$ 
be the locus over which $\mu_{\rhobar}(\sigma)$ is strictly positive.
Understanding the summation in the Breuil-Mézard conjecture (see 
Eq.~\eqref{eq:BM}) then amounts to understanding the set
$\DD(\ttt,\rhobar) = \DD(\ttt) \cap \DD(\rhobar)$.

It turns out that $\DD(\ttt)$ and $\DD(\rhobar)$ admit very
explicit combinatorial descriptions in terms of the parameters
$h$, $\gamma$ and $\gamma'$ we introduced earlier.
These descriptions first appeared in~\cite{BDJ,BP} and were then
simplified in~\cite{Da}.
Very roughly, once $\ttt$ (resp. $\rhobar$) is fixed, the weights
in $\DD(\ttt)$ (resp. in $\DD(\rhobar)$) are parametrized by
tuples $\underline \varepsilon =
\{\varepsilon_0, \ldots, \varepsilon_{f-1}\} \in \{0, 1\}^f$.
Each such tuple produces a Serre weight by a simple recipe and 
the set $\DD(\ttt)$ (resp. $\DD(\rhobar)$) is finally obtained
by putting together all such weights.
We underline that $\DD(\ttt)$ and $\DD(\rhobar)$ have usually
cardinality strictly less than $2^f$ because some $\underline 
\varepsilon$ may actually fail to produce a weight and it also
happens that two different $\underline \varepsilon$ lead to the 
same weight.

\paragraph*{The gene}

Building on the previous results, we gave in~\cite{CDM3} a purely 
combinatorial description of the intersection $\DD(\ttt,\rhobar)=\DD(\ttt) \cap 
\DD(\rhobar)$ in terms of the parameters $h$, $\gamma$ and $\gamma'$,
assuming that the tame inertial type is reducible. 
More precisely, setting $q = p^f$ for simplicity, we considered 
the quantity $h - (q{+}1) \gamma' \text{ mod } q^2 - 1$ and wrote 
its decomposition in radix $p$:
\begin{equation}
\label{eq:vi}
h - (q{+}1)\gamma' \equiv 
p^{2f-1} v_0 + p^{2f-2} v_1 + \cdots + p v_{2f-2} + v_{2f-1}
\pmod{q^2 - 1}.
\end{equation}
From the $v_i$'s, we then formed a periodic sequence 
$\bX = (X_i)_{i\in \Z}$ of period $2f$ assuming values in the 
finite set $\{\gA, \gB, \gAB, \gO\}$. The sequence $\bX$ is called
the \emph{gene} of $(\ttt, \rhobar)$ and it satisfies the following
rules (see \cite[Lemma~B.1.3]{CDM3}):
\begin{myitemize}
\item if $v_i = 0$ and $X_{i+1} = \gO$, then $X_i = \gAB$;
\item if $v_i = 0$ and $X_{i+1} \neq \gO$, then $X_i = \gA$;
\item if $v_i = 1$ and $X_{i+1} = \gO$, then $X_i = \gO$;
\item if $v_i = 1$ and $X_{i+1} \neq \gO$, then $X_i = \gB$;
\item if $v_i \geq 2$, then $X_i = \gO$.
\end{myitemize}
To each gene $\bX$, we attached a set $\WW(\bX)$ of \emph{combinatorial 
weights}, which are sequences of length $f$ with values in $\{0, 1\}$. 
We then proved (see \cite[Theorem 3.1.2]{CDM3}) that, if $\bX$ denotes
the gene of $(\ttt,\rhobar)$, there is a canonical bijection:
\begin{equation}
\label{eq:WX}
\WW(\bX) \stackrel\sim\longrightarrow \DD(\ttt,\rhobar).
\end{equation}
Beyond yielding an explicit description of $\DD(\ttt,\rhobar)$
and opening concrete and algorithmical perspectives on the 
Breuil-Mézard conjecture,
the above result raises new questions because it somehow shows 
that the dependence of $\DD(\ttt,\rhobar)$ in $\ttt$,
$\rhobar$ and even in the underlying prime number~$p$ itself, is 
very weak, given that the gene only retains few information about
these data.
In some sense, one can interpret the gene as the ``skeleton'' of
the pair $(\ttt, \rhobar)$ that captures its most fundamental
combinatorial properties in view of the Breuil-Mézard conjecture.
In this perspective, the construction $\bX \mapsto \WW(\bX)$ 
should be thought of as the core factory of Serre weights, while the 
bijections~\eqref{eq:WX}, for varying $\ttt$, $\rhobar$ and~$p$, 
appear as many tangible incarnations of this manufacture.

\paragraph*{Higher dimension and group-theoretic formulation}

When we are moving to higher dimensions, the numerical descriptions we 
used previously cannot continue to be that simple but, interestingly, 
they have analogues which can be formulated in the language of group 
theory. In what follows, we continue to assume that $K = \Q_{p^f}$ 
for some positive integer~$f$. 
In this setting, the relevant algebraic group is:
$$G = \big(\Res_{\oK/\Z_p}\GL_{n/\oK}\big) \times_{\Zp} \Zpbar
\simeq \prod_{\JJ}\GL_{n/\Zpbar}$$
where $\JJ$ is the set of embeddings $\oK \hookrightarrow \Zpbar$ 
(or, equivalenty, $k_K \hookrightarrow \Fpbar$) and will be 
identified with $\Z/f\Z$ in what follows.
We denote by $T$ the diagonal maximal torus of $G$. We let $R$ be the 
set of roots of $(G, T)$ and $W$ be the corresponding Weyl group. The 
Borel $B\subset G$ of upper triangular matrices determines a subset
$R^+ \subset R$ of positive roots.
The group of characters $X^\star(T)$ can be canonically identified 
with $(\Z^n)^\JJ = (\Z^n)^f$ and the Weyl group $W$ is isomorphic to
$\mathfrak S_n^f$. We shall also need the extended Weyl group
$\widetilde{W}$ of $G$, defined by $\widetilde{W} = W \ltimes X^\star(T) 
\simeq (\mathfrak S_n \ltimes \Z^n)^f$.

In this setting, a Serre weight is an (isomorphism class of) irreducible $\Fpbar$-representation 
of $\GL_n(\F_{p^f})$. It follows from a somehow classical argument of 
the theory of representation of reductive groups that Serre weights can 
be parametrized by certain characters of $G$.
Precisely, after~\cite[Lemma~9.2.4]{GHS}, we know that if we set:
\begin{align*}
X_0(T)
& = \big\{\,
      \lambda\in X^\star(T) \text{ s.t. }
      \langle\lambda,\alpha^\vee\rangle = 0 
      \text{ for all } \alpha\in R^+
     \,\big\} \\
X_1(T)
& = \big\{\,
      \lambda\in X^\star(T) \text{ s.t. }
      0\leq \langle\lambda,\alpha^\vee \rangle\leq p-1 
      \text{ for all } \alpha\in R^+
     \,\big\}
\end{align*}
and let $\pi$ denote the shift $(x_i)_{i \in \Z/f\Z} \mapsto 
(x_{i+1})_{i \in \Z/f\Z}$ on $X^\star(T) \simeq (\Z^n)^f$, there 
is a bijection:
\begin{equation}
\label{equaSerre}
F: X_1(T)/(p{-}\pi)X_0(T) \stackrel\sim\longrightarrow \DD
\end{equation}
taking a character $\lambda$ to the restriction to 
$\Res_{\oK/\Z_p}\GL_{n/\oK}(\F_p) = \GL_n(\F_{p^f})$ of the 
representation of $G(\F_p)$ induced by the algebraic representation of 
$G$ of highest weight $\lambda$.

We also have a description of tame inertial types
in terms of elements of the group~$\widetilde{W}$ \cite{LMS, LLHLM1,
BHHMS}.
Let $(s,\mu)\in \widetilde W = W \ltimes X^\star(T)$ and write
$s = (s_0,\ldots, s_{f-1})$ with $s_i \in \mathfrak S_n$ for all $i$.
Let $r$ be the order of $s_0s_{f-1}\cdots s_1\in \mathfrak S_n$. 
We consider the unramified extension $K'$ of $K$ of degree
$r$ and let $k_{K'}$ denote its residue field. We let
$\tilde{\omega}_{rf}:I_K=I_{K'}\to \Zpbar^\times$ be the Teichmüller 
lift of the Serre fundamental character of level $rf$. For $j\in\JJ$, 
we put:
$$\eta_j=\left\{
\begin{array}{ll} 
(n-1,\ldots,1,0) & j\mbox{-th coordinate}\\
(0,\ldots,0)&\mbox{elsewhere}
\end{array}\right.$$
and $\eta=\sum_{j\in \JJ}\eta_j$. We define
$\alpha'_{(s,\mu)} \in
X^\star(T)^{\Hom_{k_K\text{-alg}}(k_{K'},\Fpbar)}
\simeq X^\star(T)^r \simeq (\Z^n)^{rf}$ by:
$$\alpha'_{(s,\mu),j}=s_1^{-1}s_2^{-1}\cdots s_j^{-1}(\mu_j+\eta_j)$$
and finally set:
\begin{equation}
\label{equatau}
\tau(s,\mu+\eta)=\bigoplus_{1\leq i\leq n}\tilde{\omega}_{rf}^{\sum_{j'=0}^{rf-1} \alpha'_{(s,\mu),j',i}p^{j'}}.
\end{equation}
It is our tame inertial type.

Irreducible representations $\rhobar : G_K \to \GL_n(\Fpbar)$ 
can be encoded in a similar fashion. 
Moreover, when $n = 2$, it turns out that the explicit descriptions 
of $\DD(\ttt)$ and $\DD(\rhobar)$ we have mentioned earlier can be
rephrased in the language of group theory which was briefly sketched
above (at least under sufficiently generic assumptions).
We do not reproduce here the corresponding recipes 
(which involve the so-called \emph{$p$-dot product}) but refer to 
Propositions~2.4.2 and~2.4.3 of \cite{BHHMS} for more details.
So far, we do not have any candidate for being a plausible replacement 
of the gene when $n > 2$. However, all the above constructions tend to
show that, even though they now need to be formulated in the language 
of group theory, combinatorics is still here (and is maybe even more
ubiquitous) in higher dimensions.

\subsection{Explicit computations of $R^{\lambda,\ttt}_{\rhobar}$}
\label{ssec:explicitdef}

The Breuil-Mézard conjecture is concerned with the special fibre of 
$R^{\lambda,\ttt}_{\rhobar}$ but, of course, obtaining a complete 
description of the ring $R^{\lambda,\ttt}_{\rhobar}$ has also interest
for its own. For example, explicit presentations of some 
$R^{\lambda,\ttt}_{\rhobar}$ have been used 
by Emerton, Gee and Savitt~\cite{EGS} to prove important conjectures 
stated by Breuil in \cite{Br} about lattices in the cohomology of 
Shimura curves.
In this subsection, we outline the standard strategy that is used 
to approach $R^{\lambda,\ttt}_{\rhobar}$ and report on the results
of some explicit computations.

\paragraph*{Review on Kisin's construction of
$R^{\lambda,\ttt}_{\rhobar}$}

The main theoretical ingredient for studying deformations of 
$\rhobar$ with prescribed Hodge type and inertial type is the theory
of Breuil-Kisin, which provides a description of those deformations
by means of semi-linear algebra.
In our setting and assuming in addition that $\ttt$ is tame of
level $f$ and $K = \Q_{p^f}$ as we already did previously,
a Breuil-Kisin module is a 
projective module $\MK$ over $\Zpbar \otimes_{\Zp} \oK[[u]] \simeq
\Zpbar[[u]]^\JJ$ equipped with two additional structures:
\begin{myitemize}
\item a \emph{Frobenius map} $\varphi_\MK : \MK \to \MK$ which is
semi-linear (with respect to the endomorphism of $\Zpbar \otimes_{\Zp} 
\oK[[u]]$ acting by the identity on $\Zpbar$, by the Frobenius of $\oK$ 
and taking $u$ to $u^p$) and satisfy additional properties,
\item a \emph{descent data}, that is a linear action of the group 
$\Gal(K[\sqrt[e] p]/K) \simeq \Z/e\Z$ (with $e = p^f - 1)$ that
commutes with the Frobenius.
\end{myitemize}
A famous theorem of Kisin~\cite{Ki2} indicates that those modules are in 
correspondence with $\Zpbar$-representations of $G_K$ that become 
crystalline over the extension $K[\sqrt[e] p]$.
Besides, the Hodge type (resp. the inertial type) of the latter can be 
easily read off on the form of the Frobenius map (resp. of the descent
data) on the former.
Even better, the reduction modulo $p$ of the representation associated
to a Breuil-Kisin module $\MK$ is uniquely and entirely described by
the module $\MK \otimes_{\oK[[u]]} k_K((u))$ equipped with its
additional structures.

The Breuil-Kisin theory then looks particularly well fitted for the
study of the deformations rings $R^{\lambda,\ttt}_{\rhobar}$ and it
turns out that it indeed is.
In~\cite{Ki3}, Kisin constructed a scheme
$\GR^{\lambda,\ttt}_{\rhobar}$ parametrizing the Breuil-Kisin modules
$\MK$ of Hodge type $\lambda$, inertial type $\ttt$ and having the
additional property that $\MK \otimes_{\oK[[u]]} k_K((u))$ corresponds
to the given representation~$\rhobar$.
This scheme is moreover equipped with a morphism
$\GR^{\lambda,\ttt}_{\rhobar} \to \Spec R_{\rhobar}$
whose schematic image has closure $\Spec R^{\lambda,\ttt}_{\rhobar}$.
One should be careful however that the morphism:
\begin{equation}
\label{eq:kappa}
\kappa:
\GR^{\lambda,\ttt}_{\rhobar} \to \Spec R^{\lambda,\ttt}_{\rhobar}
\end{equation}
is \emph{not} an isomorphism in general because two $p$-torsion 
Breuil-Kisin modules may correspond to the same Galois representation. 
It is however always an isomorphism in generic fibre. The special fibre
of $\GR^{\lambda,\ttt}_{\rhobar}$ is denoted by
$\GRbar^{\lambda,\ttt}_{\rhobar}$ and is called the \emph{Kisin variety}; 
in some sense, it measures the default for $\kappa$ to be an isomorphism.

\paragraph*{Examples in dimension $2$: the generic case}

The first examples of explicit calculations of certain rings
$R^{\lambda,\ttt}_{\rhobar}$ have been carried out by Breuil and Mézard
in \cite{BM1} and \cite{BM2}.
They considered the case where $\rhobar$ is $2$-dimensional and 
absolutely irreducible, $\lambda = (0,0)$ for each embedding and
$\ttt$ is tame of level~$f$.
Under some additional assumptions of genericity on $\rhobar$, they
obtained, when $\DD(\ttt, \rhobar)$ is not empty:
\begin{equation}
\label{eq:R-BM}
R^{\lambda,\ttt}_{\rhobar} \simeq
\frac{\Zpbar[[ X_i, Y_i,  i \in \JJ_{\II},  Z_j, j \in \JJ\setminus \JJ_{\II} ]]}
{\left(  X_i Y_i - p, i \in \JJ_{\II}\right)}.
\end{equation}
for a certain subset $\JJ_\II$ of $\JJ$ (which depends on $\lambda$,
$\ttt$ and $\rhobar$).
The aforementioned genericity assumptions play a quite important role
in Breuil and Mézard's argument. 
In fact, they imply that the underlying Kisin variety is reduced to
one point, which itself ensures that the morphism $\kappa$ of 
Eq.~\eqref{eq:kappa} is an isomorphism.
The computation of $R^{\lambda,\ttt}_{\rhobar}$ then directly reduces
to that of $\GR^{\lambda,\ttt}_{\rhobar}$.

Before moving to nongeneric cases, it is important to comment on 
the subset $\JJ_\II$ which appeared in Eq.~\eqref{eq:R-BM}.
The triviality of the Kisin variety indicates that the module over 
$\Fpbar \otimes_{\Fp} k_K((u)) \simeq \Fpbar((u))^\JJ$ corresponding
to $\rhobar$ contains a
unique lattice $\MK(\rhobar)$ that is a Breuil-Kisin module of type $(\lambda,\ttt)$. 
Breuil and Mézard then defined the \emph{shape} of $\MK(\rhobar)$: 
it is a finite sequence $(g_0, \ldots, g_{f-1})$ assuming values 
in the finite set $\{\I, \II\}$ which, roughly speaking, is obtained 
by looking at the form of the matrix of $\varphi_{\MK(\rhobar)}$ in bases 
diagonalizing the action of the descent data. The set $\JJ_\II$ is then 
formed by the indices $i$ for which $g_i$ is $\II$.

The shape also plays a key role on the $\GL_2$-side of the Breuil-Mézard 
conjecture: in our setting, the cardinality of $\D(\ttt,\rhobar)$ is 
$2^{\text{Card}(\JJ_\II)}$ (which is the Hilbert-Samuel multiplicity of 
the special fibre of the ring $R^{\lambda,\ttt}_{\rhobar}$ given by 
Eq.~\eqref{eq:R-BM}) and, more precisely, we can explicitely parametrize 
the weights in $\D(\ttt,\rhobar)$ by subsets of $\JJ_\II$.

\paragraph*{Examples in dimension $2$: the nongeneric case}

Nongeneric cases are more complicated because they usually correspond
to nontrivial Kisin varieties.
In~\cite{CDM2}, we computed those Kisin varieties when, as above,
$\rhobar$ is absolutely irreducible, $\lambda = (0,0)$ for each
embedding and $\ttt$ is tame of level~$f$.
We recall that, in this setting, we
have attached to the pair $(\ttt, \rhobar)$ its gene $\bX$ (see
\S \ref{ssec:encode}).
Our results show that the Kisin variety is entirely determined by
the gene. Being a little bit more precise, we showed that 
$\GRbar^{\lambda,\ttt}_{\rhobar}$ is a closed subscheme of
$(\P^1_{\Fpbar})^{\Z/f\Z}$ defined by equations of the form:
\begin{equation}
\label{eq:varK}
\lambda_i \: x_i \: y_{i+1} = \mu_i \: x_{i+1} \: y_i
\end{equation}
where $[x_i:y_i]$ denotes the projective coordinates on the $i$-th copy 
of $\P^1_{\Fpbar}$ and $\lambda_i$ and $\mu_i$ are elements of $\{0,1\}$ 
that can be read off on the gene $\bX$.

It is important to observe that the notion of shape can be extended
to the nongeneric case as well.
Indeed, each $\Fpbar$-point $x$ of $\GRbar^{\lambda,\ttt}_{\rhobar}$
corresponds, by definition, to a Breuil-Kisin module and so has a
well defined shape $g(x) = (g_0(x), \ldots, g_{f-1}(x)) \in 
\{\I,\II\}^f$ in the sense of Breuil and Mézard.
In full generality, the shape is thus no longer a unique element 
in $\{\I,\II\}^f$ but a function on the $\Fpbar$-points of the Kisin 
variety taking values in $\{\I,\II\}^f$. 
We proved moreover that this function is lower-continuous (for the 
partial ordering on the codomain defined by $\I < \II$) and thus
defined a stratification on $\GRbar^{\lambda,\ttt}_{\rhobar}$ by
locally closed subschemes. As well as the Kisin variety, the shape 
stratification is entirely determined by the gene.

We then proposed the following conjecture.

\begin{conj}
\label{conj:CDM}~
\begin{myenumerate}[(i)]
\item
The generic fibre of $R^{\lambda,\ttt}_{\rhobar}$ is determined
by the Kisin variety equipped with its shape stratification.
\item
The ring $R^{\lambda,\ttt}_{\rhobar}$ is determined by the gene.
\end{myenumerate}
\end{conj}

\noindent
Regarding the first item of the conjecture, we were actually much 
more precise and exhibited a candidate for being the generic fibre
of $R^{\lambda,\ttt}_{\rhobar}$. Besides, our candidate is rather
explicit: it is defined as the formal neighborhood of the Kisin 
variety in a certain blow-up of $(\P^1_{\Zpbar})^{\Z/f\Z}$.
We refer to \cite[\S 5.4]{CDM2} for a complete exposition of this
construction.

Conjecture~\ref{conj:CDM} is known is most cases where the Kisin 
variety is trivial. It has also been checked in~\cite{CDM1} in one
example where the Kisin variety is isomorphic to $\P^1_{\Fpbar}$; 
in this case, the deformations ring we obtained is $\Zpbar[[X, Y, 
Z]]/(XY - p^2)$.

\paragraph*{Higher dimension and group-theoretic formulation} 

In dimension $3$, tamely potentially crystalline deformation rings for 
small Hodge-Tate weights and generic Galois representations have been 
studied in \cite{LLHLM}.
The authors also obtained explicit presentations of the deformations
rings $R^{\lambda,\ttt}_{\rhobar}$ in several cases.
The equations they found are often quite similar to the one given in
Eq.~\eqref{eq:R-BM}.
The case of rank $2$ unitary group has also been considered 
in~\cite{KoMo}. In this setting, it turns out that the explicit 
computations of the corresponding $R^{\lambda,\ttt}_{\rhobar}$'s
boil down to the case of $\GL_2$ (with an additional polarization structure on the Breuil-Kisin modules).
The final equations they get are then again similar to
Eq.~\eqref{eq:R-BM}.

In all these situations, although computations are certainly much more 
difficult, it is important to underline that the basic ingredients 
are the same: we continue to have a Kisin variety, together with a 
shape stratification and we hope that those two objects strongly 
govern the final form of the deformations space.
In full generality (\emph{i.e.} for any reductive group, even not
necessarily $\GL_n$), the Kisin variety and the shape stratification
are defined using techniques coming from group theory:
the Kisin variety is a subspace of the affine Grassmannian defined by
an explicit condition, which can be formulated in terms of the Cartan 
decomposition, while the shape function $x \mapsto g(x)$ is defined by 
means of the Iwahori-Cartan decomposition (and it now takes values in 
the Iwahori-Weyl group).
After~\cite{CaLe}, we know moreover that Kisin varieties are related
with Pappas-Zhu local models~\cite{PaZh}; in this connection, the
shape stratification corresponds to the canonical stratification by 
affine Schubert varieties.

\section{The field with one element}
\label{sec:f1}

The theory of the field with one element ($\F_1$) starts with an
observation of Tits~\cite{Ti} who notices inquiring numerical
coincidences in the theory of algebraic linear groups.
The most fundamental example is given by the group $\GL_n$ itself.
Indeed, its number of points over the finite field $\F_q$ is given by:
\begin{align*}
\Card \GL_n(\F_q)
 & = (q^n - 1) \cdot (q^n - q) \cdots (q^n - q^{n-1}) \\
 & = (q-1)^n \cdot q^{n(n-1)/2}\cdot [n]_q \cdot [n{-}1]_q \cdots [1]_q
\end{align*}
where $[i]_q = 1 + q + \cdots + q^{i-1}$ is the $q$-analogue of~$i$
and letting $q$ tends to $1$, we find:
\begin{equation}
\label{eq:cardGLdFq}
\Card \GL_n(\F_q) \sim_{q \to 1} (q-1)^n \cdot n!.
\end{equation}
What is surprising is that $n!$ can be interpreted as the cardinality
of the symmetric group $\mathfrak S_n$, which is nothing but the Weyl
group of $\GL_n$. Similar results hold more generally for a large
family of groups, including the orthogonal groups, the sympletic
groups and their scalar restrictions.
After these observations, Tits asked if these numerical matchings
could have deeper roots and proposed to build a geometry
of the so-called \emph{field with one element} with the objective 
to give a systematical and geometrical understanding of all the 
combinatorial structures and constructions which appear in the 
theory of Lie groups or algebraic groups.

Tits' vision was then popularized by Soulé who came up in~\cite{So} with
a first tentative definition of affine varieties over $\F_1$. Later on,
other constructions were proposed and the subject has attracted more and
more attention for the last two decades~\cite{De1, ToVa, CCM, CC, Co, 
LPL, Lo1, BBK}; see \cite{Lo2} for a recent review on this topic. 
Nowadays, the theory of $\F_1$ is not well estalished yet.
However, several definitions for the category of schemes over $\F_1$ 
have been proposed over the years and significant progress towards Tits' 
initial dream have been realized. Besides, geometry over $\F_1$ have 
been extended to new contexts and have nowadays close interactions with 
Arakelov geometry and $p$-adic geometry.
In particular, Bambozzi, Ben-Bassat and Kremnizer introduced 
in~\cite{BBK} analytic geometry over $\F_1$; they managed notably
to construct a model of the Fargues-Fontaine curve~\cite{FaFo} in their 
theory. Up to our knowledge, this was the first connexion between the 
theory of Galois representations (incarnated here by $p$-adic Hodge 
theory) and the field with one element.

\subsection{Clues in favor of a $1$-adic Breuil-Mézard conjecture}

The main reason why we believe that a $1$-adic version of the
Breuil-Mézard conjecture is possible is that many objects involved
in the formulation and/or the resolution of this conjecture do have
natural seeds in $\F_1$-geometry.
In this subsection, we list the most important of them and comment
on their $\F_1$-aspects.

First of all, we notice that the Weyl group $W$ and its extended version 
$\widetilde{W}$, which both play a quite important role in the 
construction of Serre's weights and inertial types, do have natural 
interpretations in characteristic one:
the Weyl group is the set of $\F_1$-points (this property is Tits' 
dream, which is the main guide of the theory) while the extended Weyl 
group may be interpreted as the set of points over $\F_1(X)$
(see Eq.~\eqref{eq:GLdF1X} in Appendix~\ref{app:f1-review})
Moreover, the recipes used for constructing Serre's weights and 
inertial types from an element of $\widetilde{W}$ 
(see Eqs.~\eqref{equaSerre} and~\eqref{equatau}) are purely 
combinatorial and it is quite likely that they can be reformulated 
by means of $\F_1$-geometry.

The equations of the Kisin varieties we obtained in~\cite{CDM2}
(see Eq.~\eqref{eq:varK}) show that all of them are defined over
$\F_1$.
Similarly, the deformations spaces computed in~\cite{BM1,BM2} 
and~\cite{LLHLM} all appear as product of discs and annuli (see
particularly Eq.~\eqref{eq:R-BM}); as a consequence, they
all come throught scalar extensions from analytic spaces over 
$\F_1$ in the sense of~\cite{BBK}.
Moreover, although the construction of blow-ups and formal
neighborhoods was not addressed in~\cite{BBK}, it looks quite
plausible that the candidates for deformation spaces we introduced
in~\cite{CDM2} are defined over $\F_1$ as well.

All the above examples show furthermore a strong uniform behaviour 
with respect to~$p$.
In the language of $\F_1$-geometry, this uniformity means that a
whole family of Kisin varieties (resp. of deformation spaces) 
parametrized by~$p$ comes by scalar extension to $\F_p$ (resp. 
to $\Q_p$) 
from a \emph{unique} variety (resp. analytic variety) over $\F_1$.
This result might suggest the existence of a common denominator of the 
theory of Kisin varieties (resp. deformations spaces) which is defined 
in characteristic one and underpins some of their features we observe
over the $p$-adics. This expectation is strengthened by the fact that
Kisin varieties have a deep group-theoretic interpretation (see last
paragraph of \S \ref{ssec:explicitdef}).
The same remark is valid for the shape stratification as well: they
have good chance to be visible in characteristic one, given that
thay are closely connected to affine Schubert varieties, which are 
themselves known to be defined over $\F_1$~\cite{LPL}.

The recipe of~\cite{CDM3}, giving a combinatorial description of the 
set of common Serre's weights in terms of the corresponding gene, has 
also a strong $\F_1$-flavour. Concretely, what we expect is that:
\begin{myenumerate}[(1)]
\item the gene is a sort of $\F_1$-encoding of the pair $(\ttt,
\rhobar)$,
\item the combinatorial weights of~\cite{CDM3} are the mirror of
a notion of Serre's weights in characteristic~$1$,
\item the association
$$\text{gene} \, \mapsto \, \text{set of combinatorial weights}$$
is the $\F_1$-incarnation of the construction $(\ttt, \rhobar)
\mapsto \DD(\ttt,\rhobar)$.
\end{myenumerate}
Beyond the justifications coming from the constructions of~\cite{CDM3}, 
we underline that there are other evidences supporting that Serre's 
weights in characteristic~$1$ should have something to do with 
combinatorial weights (which are, we recall, sequences of length~$f$ 
assuming values in $\{0, 1\}$).
Indeed, mimicing the usual definition in characteristic~$p$, we
expect Serre's weights in characteristic~$1$ to be interpreted as
$\bar \F_1$-representations (whatever it means) of the group 
$\GL_2(\F_{1^f})$. But, following Tits' vision, we can write:
$$\GL_2(\F_{1^f}) = \big(\Res_{\F_{1^f}/\F_1} \GL_2\big)(\F_1)
= \text{Weyl}\big(\Res_{\F_{p^f}/\F_p} \GL_2\big) = (\Z/2\Z)^f$$
and we already see the set $\{0,1\}$ entering into the
scene. More precisely, we can define $\Sym^k \F_1^2$ as
the set $\{X^k, X Y^{k-1}, \ldots, Y^k\}$ and let $\GL_2(\F_1) =
\Z/2\Z$ act on it by letting its unique nontrivial element operate 
by swapping $X$ and $Y$ (see Appendix~\ref{app:f1-review} for a 
justification of this definition). Similarly, if 
$\underline k = (k_0, \ldots, k_{f-1})$ is a tuple of integers, we 
let $\Sym^{\underline k} \F_1^2$ be the cartesian product of the
$\Sym^{k_i} \F_1^2$ equipped with the induced action of 
$\GL_2(\F_{1^f}) = (\Z/2\Z)^f$. It is then an easy exercise
to check that $\Sym^{\underline k} \F_1^2$ is irreducible 
(\emph{i.e.} the action is transitive) if and only if $k_i \in
\{0,1\}$ for all $i$.

We underline that Conjecture~\ref{conj:CDM} is also in line with 
the above vision: roughly speaking, it stipulates that the mapping 
$(\ttt,\rhobar) \mapsto R^{\ttt}_{\rhobar}$ descends over $\F_1$.

Combining the previous observations and being quite optimistic, 
one might hope that the Pappas-Rapoport spaces~\cite{PaRa} and/or
Emerton-Gee stacks~\cite{EmGe} themselves have a 
model over $\F_1$ and that the irreducible components of the special
fibre of the latter will be related to the set of Serre's weights in 
characteristic~$1$.

\subsection{Major challenges}

We do not hide that, if possible, devising a $1$-adic Langlands 
correspondence (or a $1$-adic Breuil-Mézard conjecture) will 
definitely not be a simple task. Actually, although the geometry
over $\F_1$ has already been developed quite a lot, a lot of 
fundamental ingredients and objects of the Langlands programme
are missing in characteristic $1$.

To start with, we notice that extensions of $\F_1$, usually
referred to as $\F_{1^n}$, have been already considered by several 
authors~\cite{So,KaSm,Co} but they have never been systematically 
studied.
Moreover, in the above references, $\F_{1^n}$ is defined as the 
cyclotomic extension of $\F_1$ whose Galois group is isomorphic
to $(\Z/n\Z)^\times$, and not $\Z/n\Z$. This means in particular
that we do not have apparently a nice analogue of the Frobenius
endomorphism, which sounds annoying.
An option for fixing this issue is to work with another version
of $\F_{1^n}$ on which we impose by design the existence of a
Frobenius of order~$n$ (see Appendix~\ref{app:f1-galois} for
a first rough proposal).

Similarly, the field of $1$-adic numbers $\Q_1$ and its
extensions have not attracted much attention so far.
We mention however that Connes introduced the ring of Witt vectors 
over $\F_1$ in~\cite{Co} but we are afraid that Connes' treatement
does not perfectly fit with our perspectives since it is eventually 
related to Banach algebras over the reals, and not over the $p$-adics.
In Appendix~\ref{app:q1-galois}, using a (certainly too) naive 
definition of $\Q_1$, we start exploring its theory of finite
extensions.

On a different note, it seems that the theory of representations of
reductive groups have not been systematically studied yet.
This could sound surprising given that representations of reductive
groups over finite fields have been attracted a lot of attention for
more than fifty years and the underlying theory includes a lot of
combinatorial parts that have good chance to be ``defined'' over
$\F_1$.

\subsection{Conclusion}

Although that, clearly, many locks still need to be unlocked, we
continue to believe that the $1$-adic Langlands correspondence is
possible.
Besides, according to us, trying to develop it could be, at the same 
time, a wonderful motivation and guide for exploring the $1$-adic world 
and for inspiring the $p$-adic Langlands correspondence by separating 
the universal combinatorial structures on the one hand and the more 
usual arithmetical properties on the other hand.
We thus warmly encourage all contributions to this topic.


\appendix
\section{Remarks on Galois theory in characteristic~$1$}

In this appendix, we start exploring the Galois properties of the field
with one element~$\F_1$ and the field of $1$-adic numbers~$\Q_1$.
Our aim is not at all to elaborate a complete and coherent theory
but to share our intuition and point out some difficulties.
In \S\ref{app:f1-review}, we review briefly the usual constructions
of geometries over $\F_1$. We then successively address the Galois 
theory of $\F_1$ and $\Q_1$ in \S \ref{app:f1-galois} and 
\S \ref{app:q1-galois} respectively.

\subsection{Brief review of geometry over $\F_1$}
\label{app:f1-review}

Most of modern theories of geometry over $\F_1$ start with defining the
category $\Vect{\F_1}$ of $\F_1$-vector spaces.
After Tits' observation that $\GL_d(\F_1)$
should be isomorphic to $\Sd$, it is tempting to define a vector
space over $\F_1$ simply as a set, this cardinality corresponding to its
dimension over $\F_1$. By definition, a $\F_1$-linear morphism $V \to W$
is a set-theoretical \emph{partially defined} function $f : V \to W$.
The direct sum (resp. the tensor product) of two vector spaces is their
disjoint union (resp. their cartesian product).
Besides, with this point of view, the standard $\F_1$-vector space of
dimension $d$, namely $\F_1^d$, is represented by the set $\{1, \ldots, d\}$;
its group of automorphisms then coincides with $\Sd$, \emph{i.e.}
$\GL_d(\F_1) = \mathfrak S^d$ as expected.

After vector spaces over $\F_1$, one introduces $\F_1$-algebras: they are,
by definition, objects in commutative monoids in the category $\Vect{\F_1}$,
\emph{i.e.} sets equipped with a partially defined law of commutative monoids
which is usually denoted with the multiplicative convention.
Examples of $\F_1$-algebras include $X^\N = \{1, X, X^2, \ldots\}$ and
$X^\Z$ which should be thought of as $\F_1[X]$ and $\F_1(X)$ respectively.
Similarly, the $\F_1$-algebra $\F_1[X_1, \ldots, X_d]$ is realized by the
monoid $X_1^\N X_2^\N \cdots X_d^\N$ whose elements are monomials in $X_1,
\ldots, X_d$.

If $M$ is a $\F_1$-algebra, it makes sense to define $M$-modules: they are
sets endowed with an action of $M$. The standard free module of rank~$d$
over $M$ is $M^{\oplus d} = \{1, \ldots, d\} \times M$ where $M$ acts by
multiplication on the second coordinate. It is an easy exercise to check
that the group of $M$-linear automorphisms of $M^{\oplus d}$ is the
semi-direct product $\Sd \ltimes (M^{\text{gp}})^d$ where
$M^{\text{gp}}$ denotes the subgroup of invertible elements of $M$.
In particular, we get:
\begin{equation}
\label{eq:GLdF1X}
\GL_d\big(\F_1(X)\big) = \Sd \ltimes \Z^d
\end{equation}
and thus obtain, at least in the case of $\GL_d$, a $\F_1$-style
interpretation of the extended Weyl group.

Until this point, (almost) all theories agree on definitions but, when we
are coming to $\F_1$-schemes, points of view start to diverge.
Chronologically, the first approach is due to Deitmar~\cite{De1}
and it closely follows the classical theory of schemes: Deitmar
introduced spectra of monoids, equipped them with a topology and a
notion of sheaves and finally glued them to get $\F_1$-schemes.
Soon after, Toen and Vaquié~\cite{ToVa} developed the functorial point
of view. Starting from the category $\Vect{\F_1}$ (or, more generally,
with an abstract monoidal symmetric category $\mathcal C$), they defined
the category $\Alg{F_1}$ as we did before, introduced a notion of Zariski
covering on it and finally defined $\F_1$-schemes as sheaves on $\Alg{F_1}$
for this Grothendieck topology. In~\cite{Ve}, Vezzani proved (in a slightly
different context) that Deitmar's construction on the one hand and
Toen-Vaquié's approach on the other hand are equivalent, in the sense
that they give rise to the same category of $\F_1$-schemes. Besides,
both viewpoints include a functor of scalar extension:
$$\Sch{\F_1} \to \Sch{\Z}, \quad X \mapsto X_\Z$$
deriving from the construction $M \to \Z[M]$ at the level of monoids.
Deitmar observed that toric varieties are defined over $\F_1$ but he
also proved that there are essentialy the only ones~\cite{De2}.

Another point of view on $\F_1$-schemes, which in some sense goes
back to Soulé's original definition, was proposed in~\cite{CC}
by Connes and Consani. They suggested to define a $\F_1$-scheme as a
triple $(\tilde A, X_\Z, e_X)$ where $\tilde A$ is a $\F_1$-algebra,
$X_\Z$ is a classical scheme and $e_X : (\Spec\tilde A)_\Z \to X_\Z$
is a morphism of schemes inducing a bijection on $k$-points for any
field~$k$. In some sense, this approach separates the purely combitorial
part, which is encoded by $\tilde A$, and the geometrical part, which
is delegated to the classical theory throught the scheme $X_\Z$.
López Peña and Lorscheid~\cite{LPL} proved that this framework is more
flexible in the sense that it allows for defining a much larger panel of
varieties over $\F_1$; those include grassmannians, split reductive groups,
Schubert varieties, \emph{etc.}
Besides, in a subsequent paper, Lorscheid~\cite{Lo1} concretized Tits'
premonition by realizing the Weyl group of a split reductive group as its
set of points over $\F_1$.

Beyond the field with one element, what we need for the purpose of this paper
is the field of $1$-adic numbers. Fortunately, this question has already been
touched in the literature by several authors.
In~\cite{Co}, Connes came with a definition of Witt vectors over $\F_1$.
However, Connes' construction looks a bit disconnected to our needs as it
comes equipped with a scalar extension functor assumption assuming values
in Banach algebras over the reals, and not over the $p$-adics.
In a different direction, Bambozzi, Ben-Bassat and Kremnizer~\cite{BBK}
laid the foundations of the theory of analytic varieties over~$\F_1$.
Roughly speaking, their construction is similar to the ones we briefly
sketched above except that, instead of starting with the category
$\Vect{\F_1}$, they consider various categories of sets $X$ equipped
with a function $\norm : X \to \R^+$ whose purpose is to model the
norm map. They showed that balls and annuli of rigid geometry come
from analytic varieties over~$\F_1$; this is thus also the case for all
the deformation spaces we listed in \S \ref{sec:BM}.

\subsection{Galois theory over $\F_1$}
\label{app:f1-galois}

It is a natural expectation that $\F_1$ should have a finite extension
of degree~$n$ for all $n$ with Galois group isomorphic to $\Z/n\Z$. In
the literature~\cite{So,KaSm,Co}, this extension $\F_{1^n}$ is usually
defined as the cyclotomic extension of $\F_1$,
that is the $\F_1$-algebra represented by the monoid $X^{\Z/n\Z}$.
However, it appears that this point of view is not perfectly suited to
our purpose for at least two reasons:
\begin{myenumerate}[(1)]
\setlength{\itemsep}{1ex}
\item \emph{The Galois theory is not the expected one.}
Indeed, the group of automorphisms of the monoid $X^{\Z/n\Z}$ is
$(\Z/n\Z)^\times$ and not $\Z/n\Z$; in particular, we do not have a
distinguished Frobenius endomorphism. This issue is maybe even more
visible when we extend scalars to $\Z$ or $\Fp$.
Indeed, with the above definition, one would get:
$$\F_{1^n} \otimes_{\F_1} \F_p = \F_p[X]/(X^n - 1)$$
which is certainly \emph{not} $\F_{p^n}$ and which besides has a different
Galois group.
\item \emph{The formation of $\F_{1^n}$-points does not behave as desired.}
If $X$ is a scheme defined over $\F_{1^n}$, it seems reasonable to expect
that the set of $\F_{1^n}$-points of $X$ agrees with the set of
$\F_1$-points of $\Res_{\F_{1^n}/\F_1} X$.
For $X = \GL_d$, this results in:
$$\GL_d(\F_{1^n}) =
\text{Weyl}\big(\Res_{\F_{p^n}/\F_p} \GL_d\big) = (\Sd)^n$$
where $p$, here, denotes any auxiliary prime number. This expectation is
reinforced by the fact that the group $(\Sd)^n$ plays a quite
important role in the Breuil-Mézard conjecture as recalled in \S
\ref{sec:BM}.
However, if one lets $\F_{1^n}$ be the $\F_1$-algebra corresponding to
the monoid $X^{\Z/n\Z}$, one would obtain $\GL_d(\F_{1^n}) = \Sd \ltimes
(\Z/n\Z)^d$ (see discussion before Eq.~\eqref{eq:GLdF1X}) which is certainly
\emph{not} isomorphic to $(\Sd)^n$ since even cardinalities differ!
\end{myenumerate}

The conclusion of these observations is that, although the cyclotomic
extension of $\F_1$ is undoubtedly an interesting object, it is probably
not the $\F_{1^n}$ we need for the applications we have in mind. Moreover,
one checks that there is unfortunately no $\F_1$-algebra meeting all our
requirements.
Instead, we propose to define from scratch a theory of $\F_{1^n}$-vector
spaces as we did previously for $\F_1$, trying as much as possible to
incorporate the Frobenius action and keep its desired properties.

\begin{definit'}
\label{def:F1nev}
A \emph{$\F_{1^n}$-vector space} is a set.
Given two $\F_{1^n}$-vector spaces $V$ and $W$, a \emph{$\F_{1^n}$-linear
morphism} $f : V \to W$ is the datum of $n$ partially defined 
set-theoretical
functions $f_1, \ldots, f_n : V \to W$.
\end{definit'}

Again, the standard $\F_{1^n}$-vector space of dimension $d$ is represented
by the set $\{1, \ldots, d\}$. We denote it by $(\F_{1^n})^d$, or
simply $\F_{1^n}$ when $d = 1$, in what follows. It is obvious from the
definition that the automorphism group of $(\F_{1^n})^d$ is $(\Sd)^n$,
\emph{i.e.} we have the expected equality $\GL_d(\F_{1^n}) = (\Sd)^n$.
Similarly, extending the definition of $\F_1(X)$ to our new setting,
one checks that $\GL_d(\F_{1^n}(X)) = (\Sd \ltimes \Z^d)^n$.

Moreover, we have an obvious scalar extension functor
$\Vect{\F_1} \to \Vect{\F_{1^n}}$ acting on objects by $V \mapsto V$
and on morphisms by $f \mapsto (f, f, \ldots, f)$. In what follows we
shall use the notation $\F_{1^n} \otimes_{\F_1} V$ to denote the scalar
extension of $V$ from $\F_1$ to $\F_{1^n}$.
Regarding scalar restriction, there are two different options to define
it, namely:
$$\begin{array}{rcl}
\aRes_{\F_{1^n}/\F_1} : \,\, \Vect{\F_{1^n}} & \longrightarrow & \Vect{\F_1} \smallskip \\
V & \mapsto & V^{\oplus n} \medskip \\
\mRes_{\F_{1^n}/\F_1} : \,\, \Vect{\F_{1^n}} & \longrightarrow & \Vect{\F_1} \smallskip \\
V & \mapsto & V^{\otimes n}.
\end{array}$$
(We recall that the direct sum and the tensor product over $\F_1$ are
defined as the disjoint union and the cartesian product respectively.)
The functor $\aRes_{\F_{1^n}/\F_1}$ (resp. $\mRes_{\F_{1^n}/\F_1}$) will be referred to as the \emph{additive}
(resp. the \emph{multiplicative}) scalar restriction from $\F_{1^n}$ to
$\F_1$; hence the notation. Both versions look interesting given than
they both appear as adjoints of the scalar extensions. Precisely, for
$V \in \Vect{\F_{1^n}}$ and $W \in \Vect{\F_1}$, we have:
\begin{align*}
\Hom_{\Vect{\F_{1^n}}}(V,\,\F_{1^n} \otimes_{\F_1} W)
 & = \Hom_{\Vect{F_1}}(\aRes_{\F_{1^n}/\F_1}(V),\,W), \\
\Hom_{\Vect{\F_{1^n}}}(\F_{1^n} \otimes_{\F_1} W,\,V)
 & = \Hom_{\Vect{\F_1}}(W,\,\mRes_{\F_{1^n}/\F_1}(V)).
\end{align*}
Besides, $\aRes_{\F_{1^n}/\F_1}(V)$ and $\mRes_{\F_{1^n}/\F_1}(V)$ are 
both equipped with a Frobenius which
acts by permuting cyclically the summands/factors. More concretely,
the Frobenius action on $\aRes_{\F_{1^n}/\F_1}(V) \simeq \Z/n\Z \times V$ is given by
$(i,x) \mapsto (i{+}1, x)$ and it is given on $\mRes_{\F_{1^n}/\F_1}(V) \simeq V^n$ by
$(x_1, \ldots, x_n) \mapsto (x_2, \ldots, x_n, x_1)$.
We observe in particular that $\aRes_{\F_{1^n}/\F_1}(\F_{1^n}) \simeq \Z/n\Z$.
In this sense, our definition meets the most standard presentation of
$\F_{1^n}$~\cite{So,KaSm,Co}; the main difference is that we do not
retain the group structure on $\Z/n\Z$ but replace it by a Frobenius
structure given by the shift.
This slight modification in the point of view is actually enough to
retrieve a cyclic Galois group of order~$n$.

\begin{prop'}
\label{prop:F1n-galois}
The group of automorphisms of $\aRes_{\F_{1^n}/\F_1}(\F_{1^n})$ commuting with the
Frobenius action is the cyclic group of order~$n$ generated by the
Frobenius.
\end{prop'}

\begin{proof}
An automorphism of $\aRes_{\F_{1^n}/\F_1}(\F_{1^n})$ is, by definition, a bijection
$f : \Z/n\Z \to \Z/n\Z$. Requiring that it commutes with the Frobenius
amounts to saying that $f(x+1) = f(x) + 1$ for all $x \in \Z/n\Z$.
Clearly, a function satisfying this condition must be of the form
$x \mapsto x + a$, \emph{i.e.} $f$ is a power of the Frobenius.
\end{proof}

Proposition~\ref{prop:F1n-galois} does not extend \emph{verbatim} if we 
replace $\aRes_{\F_{1^n}/\F_1}$ by $\mRes_{\F_{1^n}/\F_1}$; indeed, given that $\mRes_{\F_{1^n}/\F_1}(\F_{1^n})$ is reduced to 
one point, its group of automorphisms is trivial as well. However, we 
can recover the expected group if we consider all objects $V \in 
\Vect{\F_{1^n}}$ at the same time: the group of Frobenius-preserving 
automorphisms of the \emph{functor} $\mRes_{\F_{1^n}/\F_1}$ is cyclic of order~$n$ and
generated by the Frobenius.

Passing to the limit, we can similarly set up a theory of vector
spaces over $\bar \F_1 = \varinjlim_n \F_{1^n}$.

\begin{definit'}
A \emph{$\bar \F_1$-vector space} is a set.
Given two $\bar \F_1$-vector spaces $V$ and $W$, a \emph{$\bar \F_1$-linear
morphism} $f : V \to W$ is the datum of a sequence $(f_i)_{i \geq 0}$ of
partially defined functions from $V$ to $W$ such that for all $x \in V$,
the sequence $(f_i(x))_{i \geq 0}$ is periodic.
\end{definit'}

As before, we have a scalar extension functor $\Vect{\F_1} \to \Vect{\bar 
\F_1}$ which acts trivially on objects and takes a morphism $f$ in 
$\Vect{\F_1}$ to the constant sequence $(f, f, \ldots)$. More generally, 
there is a functor $\Vect{\F_{1^n}} \to \Vect{\bar \F_1}$ mapping a 
$\F_{1^n}$-linear morphism $(f_1, \ldots, f_n)$ to the sequence 
$(f_{i \text{ mod } n})_{i \geq 0}$.
If $V$ is finite dimensional over $\bar \F_1$ (\emph{i.e.} if $V$ 
is a finite set), any $\bar \F_1$-linear morphism with domain $V$
comes from a $\F_{1^n}$-linear morphism for some~$n$.
Restriction of scalars also exist in this context. 
Writing $\hat\Z = \varprojlim_n \Z/n\Z$, they are given by:
$$\begin{array}{rcl}
\aRes_{\bar \F_1/\F_1} : \quad \Vect{\bar \F_1} & \longrightarrow & \Vect{\F_1} \\
V & \mapsto & \hat\Z \times V \medskip \\
\mRes_{\bar \F_1/\F_1} : \quad \Vect{\bar \F_1} & \longrightarrow & \Vect{\F_1} \\
V & \mapsto & \{\,\text{periodic sequences with values in } V\,\}.
\end{array}$$
Furthermore, $\aRes_{\bar \F_1/\F_1}(V)$ and $\mRes_{\bar \F_1/\F_1}(V)$ are both equipped with a 
Frobenius endomorphism: on $\aRes_{\bar \F_1/\F_1}(V)$, it is $(i,x) \mapsto
(i{+}1, x)$ while it acts by shifting the sequence by $1$ on 
$\mRes_{\bar \F_1/\F_1}(V)$. One checks that $\Aut_\varphi(\aRes_{\bar \F_1/\F_1}(\bar \F_1)) \simeq
\Aut_\varphi(\aRes_{\bar \F_1/\F_1}) \simeq \Aut_\varphi(\mRes_{\bar \F_1/\F_1}) \simeq \hat\Z$ where
$\Aut_\varphi$ means the Frobenius-preserving automorphisms.

\subsection{Galois theory over $\Q_1$}
\label{app:q1-galois}

In what follows, we view $\Z_1$ (resp. $\Q_1)$ as the 
$\F_1$-analytic algebra (in the sense of~\cite{BBK}) corresponding
to the monoid $\varpi^\N$ (resp. $\varpi^\Z$) endowed with the norm 
$\Vert \varpi^v \Vert = r^v$ where $r$ is a fixed real number in $(0, 1)$. 
Here $\varpi$ is a formal notation for the uniformizer of $\Z_1$ and 
does not have further meaning.
Our definition of $\Z_1$ might sound too naive as it seems 
to identify $\Z_1$ with $\F_1[[\varpi]]$; however, at least for the 
properties we want to illustrate in this appendix, making this 
confusion will not have undesirable consequences.

We now aim at defining several families of extensions of $\Q_1$
and studying their Galois properties.
We start with unramified extensions: we set $\Q_{1^n} = \F_{1^n} 
\otimes_{\F_1} \Q_1$ for any positive integer $n$ and
$\Q_1^\ur = \bar \F_1 \otimes_{\F_1} \Q_1$. They are equipped
with a Frobenius structure coming from the Frobenius on $\F_{1^n}$
(resp. $\bar \F_1$) and with a structure of $\Q_1$-algebra 
materialized by the multiplication morphism by $\varpi$ at the
level of monoids.
One checks that the group of automorphisms of $\aRes_{\F_{1^n}/\F_1} 
(\Q_{1^n})$ (resp. of $\aRes_{\bar \F_1/\F_1} (\Q_1^\ur)$) commuting 
with both structures is isomorphic to $\Z/n\Z$ (resp. to $\hat\Z$).
We then get the expected Galois group.

We now come to the analogue of the tower of tamely ramified
extensions. When $p$ is an actual prime number, this tower is
obtained by extracting $e$-th roots of the uniformizer for $e$
coprime with~$p$. There exists an obvious analogue of this 
construction over $\Q_1$: for any positive integer $e$ (without
any condition of coprimality), we consider the $\F_1$-analytic
algebra $\Q_1[\sqrt[e]{\varpi}]$ defined by the underlying 
monoid $\varpi^{(1\!/\!e) \cdot \Z}$ equipped with the norm 
$\Vert \varpi^v \Vert = r^v$ ($v \in \frac 1 e \Z$). Clearly
$\Q_1[\sqrt[e]{\varpi}]$ is an extension of $\Q_1$ and we can
define generally 
$\Q_{1^n}[\sqrt[e]{\varpi}] = \F_{1^n} \otimes_{\F_1} \Q_1[\sqrt[e]{\varpi}]$
and
$\Q_1^\ur[\sqrt[e]{\varpi}] = \bar \F_1 \otimes_{\F_1} \Q_1[\sqrt[e]{\varpi}]$.
It is also possible to take the limit on $e$ and define 
$\Q_1[\sqrt[\infty]{\varpi}]$ as the $\F_1$-analytic algebra
associated to the monoid $\varpi^\Q$ with norm $\Vert \varpi^v 
\Vert = r^v$.
We write $\Q_1^\tr = \bar \F_1 \otimes_{\F_1} \Q_1[\sqrt[\infty]
{\varpi}]$; it is our candidate for being the maximal tamely
ramified extension (or even an algebraic closure?) of $\Q_1$.

Devising a decent Galois theory for the extension $\Q_1^\tr /
\Q_1^\ur$ looks more difficult. Indeed, given that a morphism
of monoids $\Q \to \Q$ which acts by the identity on $\Z$ needs
to be trivial, we conclude that there is no nontrivial morphism
of $\Q_1^\ur$-algebras of $\Q_1^\tr$ in the sense of
Definition~\ref{def:F1nev}.
In order to explain how this issue can be tackled, it will be more 
convenient to work with finite extensions. For any positive integer
$n$, we define $K_n = \F_{1^n} [\sqrt[n]{\varpi}]$; it is the 
$\F_{1^n}$-algebra represented by the monoid $\eta^\Z$ where we have
set $\eta = \sqrt[n]{\varpi}$ for simplicity.
As we said earlier, there are no nontrivial automorphisms of
$\Q_{1^n}$-algebras of $K_n$. The subtlety is that such automorphisms
do exist after restricting scalars to $\F_1$. An explicit example is
given by the morphism
$$\sigma_n : \aRes_{\F_{1^n}/\F_1}(K_n) \to
\aRes_{\F_{1^n}/\F_1}(K_n)$$
corresponding to the map:
$$\begin{array}{rcl}
\sigma_n^\sharp: \,\,
\Z/n\Z \times \eta^\Z & \longrightarrow & \Z/n\Z \times \eta^\Z 
\smallskip \\
(i,\,\eta^j) & \mapsto & (i{+}j,\,\eta^j).
\end{array}$$
One checks that $\sigma_n^\sharp$ is a morphism of monoids (where 
the factor $\Z/n\Z$ is endowed with its additive structure) acting
by the identity on  $\aRes_{\F_{1^n}/\F_1}(\Q_{1^n})$.
Therefore, $\sigma_n$ has all the virtues to be considered as an 
element of the
Galois group $\Gal(K_n/\Q_{1^n})$, although it is not clear to us,
here, why we need to retain the monoid structure on $\Z/n\Z$
whereas we argued earlier that it should be discarded in favor of
the Frobenius structure.
In any case, we have the following proposition.

\begin{prop'}
\label{prop:Galtr}
The group of automorphisms of monoids of $\Z/n\Z \times \eta^\Z$
acting trivially on the submonoid $\Z/n\Z \times \varpi^\Z$ is
cyclic of order $n$, generated by $\sigma_n^\sharp$.
\end{prop'}

\begin{proof}
Let $f$ be an automorphism of $\Z/n\Z \times \eta^\Z$ satisfying
the conditions of the proposition. By assumption $f$ fixes $(1,1)$
and $(0,\varpi)$. Write $f(0,\eta) = (a, \eta^b)$ with $a \in
\Z/n\Z$ and $b \in \Z$. Since $f$ is a morphism of monoids, we
must have $(na, \eta^{nb}) = (0,\varpi)$, showing that $b = 1$.
Hence $f$ takes the form $(i,\eta^j) \mapsto (i{+}aj,\, \eta^j)$
and the proposition follows.
\end{proof}

After what precedes, we are tempted to write:
$$\Gal(K_n/\Q_{1^n}) = \langle \sigma_n \rangle \simeq \Z/n\Z.$$
Moreover, noticing that $\sigma_n$ commutes with the Frobenius, we
conclude that:
$$\Gal(K_n/\Q_1) = \langle \varphi, \sigma_n \rangle \simeq
(\Z/n\Z)^2.$$
Passing to the limit, we would end up with 
$\Gal(\Q_1^\tr/\Q_1^\ur) 
\simeq \hat\Z$ and $\Gal(\Q_1^\tr/\Q_1) \simeq \hat \Z^2$.

In order to give more credit to this conclusion, we would like to make 
the comparison with the classical case of $\Q_p$ (where $p$ is an actual 
prime number). For each positive integer~$n$, we set $K_{p,n} = 
\Q_{p^n}[p^{1/(p^n-1)}]$; it is the maximal totally and tamely ramified 
Galois extension of $\Q_{p^n}$. As a consequence, the extensions 
$K_{p,n}$ are cofinal in the maximal tamely ramified extension of 
$\Q_p$. Besides, the Galois group of $K_{p,n}/\Q_p$ sits in the 
following exact sequence:

\medskip

\hfill%
\begin{tikzpicture}[xscale=3, yscale=1]
\node (A) at (0.3,0) { \ph $1$ };
\node (B) at (1,0) { \ph $\Gal(K_{p,n}/\Q_{p^n})$ };
\node (C) at (2.03,0) { \ph $\Gal(K_{p,n}/\Q_p)$ };
\node (D) at (3,0) { \ph $\Gal(\Q_{p^n}/\Q_p)$ };
\node (E) at (3.7,0) { \ph $1$ };
\draw[-latex] (A)--(B);
\draw[-latex] (B)--(C);
\draw[-latex] (C)--(D);
\draw[-latex] (D)--(E);
\node[scale=0.9] at (1,-0.9) { \ph cyclic of };
\node[scale=0.9] at (1,-1.3) { \ph order $p^n{-}1$ };
\draw[->] (1,-0.6)--(B);
\node[scale=0.9] at (3,-0.9) { \ph cyclic of };
\node[scale=0.9] at (3,-1.3) { \ph order $n$ };
\draw[->] (3,-0.6)--(D);
\end{tikzpicture}%
\hfill\null

\smallskip

\noindent
which admits a section and provides a presentation of 
$\Gal(K_{p,n}/\Q_p)$ as a semi-direct product $\Z/n\Z \ltimes
\Z/(p^n{-}1)\Z$ where $a \in \Z/n\Z$ acts on $\Z/(p^n{-}1)\Z$
by multiplication by $p^a$.

Setting naively $p=1$ in what precedes, we find $p^n{-}1 = 0$
which does not really make sense since $\Gal(K_n/\Q_{1^n})$ cannot
be decently of cardinality zero.
Remembering what we did in Eq.~\eqref{eq:cardGLdFq}, we instead
write the factorization:
$$p^n - 1 = (p-1) \cdot (1 + p + \cdots + p^{n-1}) = 
(p-1) \cdot [n]_p.$$
At the level of groups, the above factorization reflects the fact 
that $\Gal(K_{p,n}/\Q_{p^n})$ admits a subgroup of order $[n]_p$, 
namely $\Gal(K_{p,n}/K_{p,1})$.
When $p$ goes to $1$, the factor $[n]_p$ converges to $n$ and then,
passing to the limit, we expect the group $\Gal(K_n/K_1) =
\Gal(K_n/\Q_{1^n})$ to be
cyclic of order~$n$, which is exactly what we have found earlier.
Furthermore, when $p$ tends to $1$, the action of $\Z/n\Z$ on 
$\Z/(p^n{-}1)\Z$ (and consequently on all its subgroups) becomes
trivial, confirming our prediction that $\Gal(K_n/\Q_1)$ should be
a direct product of $\Gal(K_n/\Q_{1^n})$ by $\Gal(\Q_{1^n}/\Q_1)$.

\begin{rem'}
There is however one small annoying point in what we have said:
why is it legitimate to get rid of the factor $(p{-}1)$?
If instead of discarding it without further discussion, we try to
keep it, we come to the conclusion that there should be between
$\Q_{1^n}$ and $K_1$ an extension of degree $0$ or, say, of
infinitesimal degree.
This suggests that the extensions $\Q_{1^n}$ and $K_1$ need to be 
considered as different objects, which could be a way to explain 
that the prefactor $\Z/n\Z$ should be endowed with its Frobenius 
structure in the former case and with its monoid structure in the
latter one (see discussion before Proposition~\ref{prop:Galtr}).

\vspace{-\parskip}

In the similar fashion that the factor $(q{-}1)^n$ in 
Eq.~\eqref{eq:cardGLdFq} corresponds to the $n$-dimensional torus 
of $\GL_n$, it is tempting to interpret the Galois group of the ghost 
infinitesimal extension $K_1/\Q_{1^n}$ as the algebraic group $\Gm$
over $\F_1$.
Similarly, it sounds plausible to interpret the cyclic group $\Z/n\Z$
(which is supposed to be the Galois group of $K_n/K_1$) as the group of 
$\F_1$-points of an algebraic group, maybe $\aRes_{\F_{1^n}/\F_1}(\Ga)$.
All of this suggests that $\Gal(K_n/\Q_{1^n})$ could just be the
pale reflection of an algebraic group $\GGal(K_n/\Q_{1^n})$ defined
over $\F_1$ and sitting in an exact sequence of the form:
$$1 \to \aRes_{\F_{1^n}/\F_1}(\Ga) \to \GGal(K_n/\Q_{1^n})
\to \Gm \to 1.$$
And similarly, passing to the limit:
$$1 \to \aRes_{\bar \F_1/\F_1}(\Ga) \to \GGal(\Q_1^\tr/\Q_1^\ur)
\to \Gm \to 1.$$
Beyond its own interest, this interpretation would provide us with a
natural Frobenius structure (given by $i \mapsto i+1$) on 
$\Gal(\Q_1^\tr/\Q_1^\ur) \simeq \hat\Z$ after taking $\F_1$-points.
\end{rem'}

\end{document}